\documentclass{gtmon_a}
\pdfoutput=1
\usepackage[all]{xy}
\usepackage{amscd}

\proceedingstitle{The Zieschang Gedenkschrift}
\conferencestart{5 September 2007}
\conferenceend{8 September 2007}
\conferencename{Conference in honour of Heiner Zieschang}
\conferencelocation{Toulouse, France}

\editor{Michel Boileau}
\givenname{Michel}
\surname{Boileau}

\editor{Martin Scharlemann}
\givenname{Martin}
\surname{Scharlemann}

\editor{Richard Weidmann}
\givenname{Richard}
\surname{Weidmann}

\title[Minimizing coincidence numbers of maps into 
       projective spaces]{Minimizing coincidence 
       numbers\\of maps into projective spaces}

\author{Ulrich Koschorke}
\givenname{Ulrich}
\surname{Koschorke}
\address{Universit\"at Siegen\\Emmy Noether Campus\\\newline
Walter-Flex-Str. 3\\D-57068 Siegen\\Germany}
\email{koschorke@mathematik.uni-siegen.de}
\urladdr{http://www.math.uni-siegen.de/topology/}

\volumenumber{14}
\issuenumber{}
\publicationyear{2008}
\papernumber{17}
\startpage{373}
\endpage{391}

\doi{}
\MR{}
\Zbl{}



\received{21 March 2006}
\revised{}
\accepted{26 January 2007}
\published{29 April 2008}
\publishedonline{29 April 2008}
\proposed{}
\seconded{}
\corresponding{}
\version{}

\arxivreference{math.AT/0606024}

\makeatletter
\def\cnewtheorem#1[#2]#3{\newtheorem{#1}{#3}[section]
\expandafter\let\csname c@#1\endcsname\c@subsection}
\makeatother

\let\xysavmatrix\xymatrix
\def\xymatrix{\disablesubscriptcorrection\xysavmatrix}
\AtBeginDocument{\let\bar\wbar\let\tilde\wtilde\let\widetilde\wtilde\let\hat\what%
\let\overline\wbar}

\swapnumbers
\cnewtheorem{thm}[subsection]{Theorem}
\cnewtheorem{prop}[subsection]{Proposition}
\cnewtheorem{cor}[subsection]{Corollary}

\theoremstyle{definition}
\cnewtheorem{definition}[subsection]{Definition}
\cnewtheorem{example}[subsection]{Example}
\cnewtheorem{tble}[subsection]{Table}

\cnewtheorem{quest}[subsection]{Question}
\cnewtheorem{notcon}[subsection]{Notation and conventions}
\cnewtheorem{remark}[subsection]{Remark}
\cnewtheorem{numberonly}[subsection]{}
\makeatletter
  \let\c@equation\c@subsection
  \def\theequation{\thesubsection}
\makeatother
\makeautorefname{tble}{Table}

\newcommand{\pr}{\operatorname{pr}}
\newcommand{\incl}{\operatorname{incl}}

\newcommand{\im}{\operatorname{im}}

\newcommand{\scirc}{{\scriptscriptstyle{\circ}}}

\begin{document}

\begin{htmlabstract}
In this paper we continue to study (&lsquo;strong&rsquo;) Nielsen coincidence
numbers (which were introduced recently for pairs of maps between
manifolds of arbitrary dimensions) and the corresponding minimum
numbers of coincidence points and pathcomponents. We explore
compatibilities with fibrations and, more specifically, with covering
maps, paying special attention to selfcoincidence questions. As a
sample application we calculate each of these numbers for all maps
from spheres to (real, complex, or quaternionic) projective
spaces. Our results turn out to be intimately related to recent work
of D Gon&ccedil;alves and D Randall concerning maps which can be deformed
away from themselves but not by small deformations; in particular,
there are close connections to the Strong Kervaire Invariant One
Problem.
\end{htmlabstract}

\begin{abstract}
In this paper we continue to study (`strong') Nielsen coincidence
numbers (which were introduced recently for pairs of maps between
manifolds of arbitrary dimensions) and the corresponding minimum
numbers of coincidence points and pathcomponents. We explore
compatibilities with fibrations and, more specifically, with covering
maps, paying special attention to selfcoincidence questions. As a
sample application we calculate each of these numbers for all maps
from spheres to (real, complex, or quaternionic) projective
spaces. Our results turn out to be intimately related to recent work
of D Gon\c calves and D Randall concerning maps which can be deformed
away from themselves but not by small deformations; in particular,
there are close connections to the Strong Kervaire Invariant One
Problem.
\end{abstract}

\begin{asciiabstract}
In this paper we continue to study (`strong') Nielsen coincidence
numbers (which were introduced recently for pairs of maps between
manifolds of arbitrary dimensions) and the corresponding minimum
numbers of coincidence points and pathcomponents. We explore
compatibilities with fibrations and, more specifically, with covering
maps, paying special attention to selfcoincidence questions. As a
sample application we calculate each of these numbers for all maps
from spheres to (real, complex, or quaternionic) projective
spaces. Our results turn out to be intimately related to recent work
of D Goncalves and D Randall concerning maps which can be deformed
away from themselves but not by small deformations; in particular,
there are close connections to the Strong Kervaire Invariant One
Problem.
\end{asciiabstract}

\maketitle

\section{Introduction and statement of results}\label{sec1} 

This work is dedicated to the memory of Heiner Zieschang who
contributed so substantially to many areas of topology and, in
particular, to coincidence theory in codimensions~$0$ (see for example
Bogatyi, Gon{\c{c}}alves, Kudryavtseva and Zieschang \cite{BGZ2,GKZ}
and the references given there, as well as their very clear and helpful
survey article \cite{BGZ1}).

In this paper we study coincidences in higher codimensions.  As in our
paper \cite{K6} we will use geometric methods which involve
(nonstabilized) normal bordism theory and path spaces.  It was shown
in \cite{K4} that a similar approach can also be applied fruitfully to
certain knotting and linking phenomena -- another subject which Heiner
Zieschang investigated for many years.

Let $f_1, f_2 \co M \to N$ be (continuous) maps between smooth connected manifolds without boundary of arbitrary
positive dimensions $m$ and $n$, $M$ being compact. Our interest centers around coincidence spaces such as
\begin{equation}\label{1.1}
C (f_1, f_2) = \{ x \in M \mid f_1 (x) = f_2 (x)\} .
\end{equation}           
We want to determine the {\bf m}inimum number of {\bf c}oincidence points
\begin{equation}\label{1.2}
MC (f_1, f_2) \ := \ \min \{\# C (f'_1, f'_2)\mid f'_1 \sim  f_1, f'_2 \sim f_2\}
\end{equation}            
and the {\bf m}inimum number of {\bf c}oincidence {\bf c}omponents
\begin{equation}\label{1.3}
MCC (f_1, f_2) :=\ \min \{ \# \pi_0 (C (f'_1, f'_2)) \mid f'_1 \sim f_1, f'_2 \sim f_2\}
\end{equation}             
(compare \cite{K6} for the details concerning these and the following definitions).

Generically the map $(f_1, f_2) \co M \to N \times N$ is smooth and transverse to the diagonal $\Delta = \{ (y, y) \in N \times N \mid y \in N\}$. Then the coincidence locus
\begin{equation}\label{1.4}
C = C (f_1, f_2) \ = \ (f_1, f_2)^{- 1} (\Delta) \ \  \subset \ M
\end{equation}              
is a closed smooth $(m - n)$--dimensional submanifold of $M$. Its normal bundle $\nu (C, M)$ is described by a canonically arising vector bundle isomorphism
\begin{equation}\label{1.5}
\overline g^\# \co \nu (C, M) \cong f_1^* (TN) | C \ .
\end{equation}                
Moreover the obvious fiber projection from the space
\begin{equation}\label{1.6}
E (f_1, f_2) = \{ (x, \theta) \in M \times P (N) \mid \theta (0) = f_1 (x), \theta (1) = f_2 (x) \}
\end{equation}                 
(where $P (N)$ denotes the set of continuous paths, equipped with the compact--open topology)) to $M$ allows a canonical section $\widetilde g$ over $C$:
\begin{equation}\label{1.7}
\widetilde g (x) := (x, \ \text{constant path at} \  f_1 (x) = f_2 (x)), \ \ \ x \in C \ .
\end{equation}                   
The resulting bordism class
\begin{equation}\label{1.8}
\omega^\# (f_1, f_2) \ = \ [C, \widetilde g, \overline g^\#] \ \in \ \Omega^\# (f_1, f_2)
\end{equation}                    
in the bordism set of such triples is our key invariant (see \cite{K6}
for details).

\begin{definition}\ \label{1.9}

(i)\qua We call a pathcomponent $A$ of the space $E (f_1, f_2)$ {\em strongly essential} if the corresponding (partial) bordism class
\begin{equation*}
\omega^\#_A (f_1, f_2) \ := \ [C_A := \widetilde g^{- 1} (A), \ \widetilde g | C_A, \ \overline g^\# | ]
\end{equation*}
is nontrivial.

(ii)\qua We define the {\em strong Nielsen number} $N^\# (f_1, f_2)$ of $f_1$ and $f_2$ to be the number of strongly essential pathcomponents $A \in \pi_0 (E (f_1, f_2))$.
\end{definition}

\begin{thm}{\rm (cf \cite[1.2 and 3.1]{K6})}\qua\label{1.10}The finite 
number $N^\# (f_1, f_2) = N^\# (f_2, f_1)$ depends only on the
homotopy classes of $f_1$ and $f_2$. We have:
\begin{equation*}
0 \ \le \ N^\# (f_1, f_2) \ \le \ MCC (f_1, f_2) \le MC (f_1, f_2) \
\le \infty \ .
\end{equation*}
If $n \ne 2$, then also $MCC (f_1, f_2) \le \# \ \pi_0 (E (f_1, f_2)$;
if $(m, n) \ne (2, 2)$, then\break $MC (f_1, f_2) \le \# \pi_0 (E (f_1,
f_2)$ or $MC (f_1, f_2) = \infty$.
\end{thm}

(Recall from \cite[2.1]{K3} that the cardinality $\# \pi_0 (E (f_1, f_2))$ can be interpreted as the cardinality of a well known Reidemeister set.)

\begin{definition}\label{1.11}
We call the pair $(f_1, f_2)$ {\em loose} if the maps $f_1, f_2$ can be deformed away from one another or, equivalently, if $M(C)C (f_1, f_2) = 0$.

We have the obvious implications
\begin{equation}\label{1.12}
(f_1, f_2) \ \text{is loose} \ \ \Longrightarrow \ \ \omega^\# (f_1, f_2) = 0 \ \ \Longrightarrow \ \ N^\# (f_1, f_2) \ = \ 0 .
\end{equation} 

\end{definition}

\begin{quest}\label{1.13}
Do the converse implications hold?  More generally: does the minimum coincidence number $MCC (f_1, f_2)$ coincide with the strong Nielsen number $N^\# (f_1, f_2)$?

For example if $N^\# (f_1, f_2) = 0$ then there are individual
nulbordisms for each $C_A$ (cf \ref{1.9}); but do they fit together to
yield a disjointly embedded nulbordism for {\em all} of $C = \ \amalg
\ C_A$?
\end{quest}

The answer to \fullref{1.13} is positive for maps between spheres.

\begin{example}(cf \cite[1.12]{K6})\label{1.14}\qua Consider maps $f_1, f_2 \co S^m \to S^n$ where $m, n \ge 1$, and let  $A$ denote the antipodal involution. Then
\begin{equation*}
MCC (f_1, f_2) = N^\# (f_1, f_2) = \left\lbrace \begin{matrix}\ 0\hfill  &\hbox{if} \ f_1 \sim A \scirc f_2 \ ; \\
\ \# \pi_0 (E (f_1, f_2)) &\hbox{otherwise} \ .\hfill \\  \end{matrix} \right.
\end{equation*}
\end{example}

\begin{remark}\label{1.15}
Clearly $N^\#, MCC$ and $MC$ coincide for all $f_1, f_2 \co M \to N$ whenever $m < n$ or, in case $M = S^m$, whenever $m = 1$ or $n = 1$ (compare eg \cite[1.3]{K6}). Thus we will be mainly interested in situations where $m \ge n \ge 2$.
\end{remark}

In this paper we give a new interpretation of our Nielsen numbers in terms of (liftings to) covering spaces (see \fullref{sec3} below). As a sample application we compute Nielsen and minimum coincidence numbers for all pairs of maps from spheres to projective spaces.

Let $\mathbb K = \mathbb R, \mathbb C$ or $\mathbb H$ denote the field of real, complex or quaternionic numbers, and let $d = 1, 2,$ or $4$ be its real dimension. Let $\mathbb K P (n')$ and $V_{n'+1, 2} (\mathbb K)$, respectively, denote the corresponding space of lines and of orthonormal 2--frames, respectively, in $\mathbb K^{n'+1}$. The real dimension of $N = \mathbb K P (n')$ is $n := d \cdot n'$.

Consider the diagram
\begin{equation}\label{1.16}
\begin{CD}
\dots \longrightarrow  \pi_m (V_{n'+1, 2} (\mathbb K)) @>{p_{\mathbb K *}}>> \pi_m (S^{n+d-1})  @>{\partial_{\mathbb K}}>> \pi_{m -1} (S^{n -1}) \longrightarrow \dots \\
@.                        @VV{p_*}V                                                                @VV{E}V         \\
@.                    \pi_m (\mathbb K P (n'))                                 @.               \pi_m (S^n)
\end{CD}
\end{equation}
determined by the canonical fibrations $p$ and $p_{\mathbb K}$; $E$ denotes the Freudenthal suspension homomorphism.

In view of \fullref{1.14} and \fullref{1.15} above (as well as of example 1.12 and the appendix in \cite{K6}) the following result determines our Nielsen and minimum numbers for all $f_1, f_2 \co S^m \to \mathbb K P (n'), \ m, n' \ge 1$. (Proofs will be given in \fullref{sec6} below).

\begin{thm}\label{1.17} 
Assume $m, n' \ge 2$. Given $[f_i] \in \pi_m (\mathbb K P (n'))$, there is a unique homotopy class $[\widetilde f_i] \in \pi_m (S^{n + d - 1})$ such that $p_* ([\widetilde f_i]) - [f_i]$ lies in the image of $\pi_m (\mathbb K P (n' - 1)), \ i = 1, 2$. (Since this image is isomorphic to $\pi_{m -1} (S^{d -1})$, we may assume that $\widetilde f_i$ is a genuine lifting of $f_i$ when $\mathbb K = \mathbb R$ or when $m > 2$ and $\mathbb K = \mathbb C)$. Define $[f'_i] := [p \scirc \widetilde f_i ] \in \pi_m (\mathbb K P (n'))$.

Each pair of homotopy classes $[f_1], [f_2] \in \pi_m (\mathbb K P
(n'))$ satisfies precisely one of the seven conditions which are
listed in \fullref{1.18} below, together with the corresponding
Nielsen and minimum numbers.
\end{thm}

{\small\renewcommand{\arraystretch}{1.5}
\begin{tabular}{l|c|c|c|} \hline
{Condition } & ${\scriptstyle N^\# (f_1, f_2)}$ & ${\scriptstyle MCC (f_1, f_2)}$ & ${\scriptstyle MC (f_1, f_2)}$ \\ \hline\hline
1) \ ${\textstyle f'_1 \, \sim \, f'_2,\ \ [\widetilde f_2] \ \in \ \ker \partial_{\mathbb K}}$  & 0 & 0 & 0 \\ \hline
2) \ ${\textstyle f'_1 \, \sim \, f'_2,\ \ [\widetilde f_2] \ \in \ \ker E \scirc \partial_{\mathbb K} -  \ker \partial_{\mathbb K}}$ & 0 & 1 & 1 \\ \hline
3) \ ${\textstyle \mathbb K \, = \, \mathbb R, \ f'_1 \sim f'_2, \ \ \widetilde f_2 \not\sim A \scirc \widetilde f_2}$ & 1 & 1 & 1 \\ \hline
4) \ ${\textstyle  \mathbb K  =  \mathbb R,} \ {\textstyle f'_1 \not\sim f'_2, \ \ [\widetilde f_1] - [\widetilde f_2]}  \in  {\textstyle E (\pi_{m -1} (S^{n -1}))} $ & 2 & 2 & 2 \\ \hline
5) \ ${\textstyle \mathbb K \, = \, \mathbb R,  \ [\widetilde f_1] - [\widetilde f_2] \ \not\in \ E (\pi_{m -1} (S^{n -1}))}$ & 2 & 2 & $\infty$ \\ \hline
6) \ ${\textstyle \mathbb K \, = \, \mathbb C \, \text{or} \, \mathbb H, \ [\widetilde f_1] \, = \, [\widetilde f_2] \ \not\in \ \ker E \scirc \partial_{\mathbb K}}$ & 1 & 1 & 1 \\ \hline
7) \ ${\textstyle \mathbb K \, = \, \mathbb C \, \text{or} \, \mathbb H, \ [\widetilde f_1] \ne [\widetilde f_2]}$ & 1 & 1 & $\infty$ \\ \hline
\end{tabular}}

\begin{tble}\label{1.18}
Nielsen and minimum coincidence numbers of all pairs of maps $f_1, f_2
\co$ $S^m \to \mathbb K P (n'), \ m, n' \ge 2$: replace each (possibly
base point free) homotopy class $[f_i]$ by a base point preserving
representative and read off the values of $N^\#$ and $M (C)C$. (Here
$f'_1 \sim f'_2$ means that $f'_1, f'_2$ are homotopic in the
basepoint free sense; $A$ denotes the antipodal map.)
\end{tble}

\begin{example}[$\mathbb K = \mathbb R, \ m = 11, \ n = 6$, compare
  Gon{\c{c}}alves and Wong {\cite[Example 2.4]{GW}}]\label{1.19}
 According to Toda \cite{T} and Paechter \cite{P} we have in \ref{1.16}
\begin{equation*}
\frac{1}{2} H \co \pi_{11} (S^6)  \overset{\cong}{\longrightarrow} \mathbb Z; \ \ \pi_{10} (S^5) \cong \mathbb Z_2; \ \  \pi_{10} (V_{7, 2}) = 0
\end{equation*}
where $H$ denotes the Hopf invariant. Thus $\partial_{\mathbb K}$ is surjective and $E \equiv 0$ here.

Given maps $f_1, f_2 \co S^{11} \to \mathbb R P (6)$, the numbers $N^\# (f_1, f_2)$ and $MCC (f_1, f_2)$ differ precisely if $f_1$ and $f_2$ are homotopic and $H (\widetilde f_i) \equiv 2 (4), \ i = 1, 2$. In this case $\omega^\# (f_1, f_2) = 0$ but the pair $(f_1, f_2)$ is not loose. \hfill $\square$
\end{example}

In \cite{K3} we studied the looseness obstruction
\begin{equation} \label{1.20}
\widetilde\omega (f_1, f_2) \ \in \ \Omega_{m -n} (E (f_1, f_2); \ \widetilde\varphi)
\end{equation}
which lies in a manageable normal bordism group (or, equivalently, in a stable homotopy group) and is often accessible to computations. However already when $N$ is a sphere $\widetilde\omega$ turned out not to be a {\em complete} looseness obstruction. In order to remedy this we introduced $\omega^\#$ as a \lq\lq desuspended\rq\rq version of $\widetilde\omega$ which captures also nonstabilized geometric coincidence data. As the counterexamples in \ref{1.19} show we have not quite desuspended far enough: if $[f_1] = [f_2] \in \pi_{11} (\mathbb R P (6))$ then the precise looseness obstruction for $(f_1, f_2)$ is the homotopy class $\partial_{\mathbb R} ([\widetilde f_1])$, but $\omega^\# (f_1, f_2)$ is only as strong as the (once!) suspended value $E (\partial_{\mathbb R} ([\widetilde f_1])$.

\fullref{1.19} implies that  the first part of the first question as well as the second question in \ref{1.13} does not always have a positive answer. However, if $M = S^m$ and $N$ is a projective space then $N^\#$ has the following property (analogous to a norm on a vector space).

\begin{thm}\label{1.21}
Given any maps $f_1, f_2 \co S^m \to \mathbb K P (n')$ $($where $m, n' \ge 1$ \ and \ $\mathbb K = \mathbb R, \mathbb C$ or $\mathbb H)$, we have:
\begin{equation*}
\omega^\# (f_1, f_2) = 0 \ \ \ \text{if and only if} \ \ \ N^\# (f_1, f_2) = 0
\end{equation*}
\end{thm}

This follows from the fact that for $m, n' \ge 2$ the Nielsen number
vanishes only in the selfcoincidence setting where $f_1 \sim f_2$ and
automatically all pathcomponents of $E (f_1, f_2)$ but one are
strongly inessential (cf \ref{5.1}).

Observe also that the weakness of the $\omega^\#$--invariant
illustrated by \fullref{1.19} occurs here only in the
selfcoincidence setting. In contrast, in the root setting (where $f_2$
is constant) $\omega^\# (f_1, *)$ is a complete looseness obstruction
for all $[f_1] \in \pi_m (\mathbb K P (n'))$ (compare \cite[6.5]{K6}).

For $\mathbb K = \mathbb R$ let us give a more systematic treatment of the case 2 in \fullref{1.18}.

\begin{thm}\label{1.22}
Given a map $f \co S^m \to \mathbb R P(n), \ m, n \ge 2$, let $\widetilde f \co S^m \to S^n$ be a lifting. Then the following conditions are equivalent:

\begin{enumerate}
\item \ $\omega^\# (f, f) = 0$,  but $(f, f)$ is not loose;
\item \ $\partial_{\mathbb R} ([\widetilde f]) \ne 0$, but $E \scirc \partial_{\mathbb R} ([\widetilde f]) = 0$;
\item \ $(\widetilde f, \widetilde f)$ is loose, but $(f, f)$ is not loose;
\item \ $MC (\widetilde f, \widetilde f) < MC (f, f)$;
\item \ $MCC (\widetilde f, \widetilde f) < MCC (f, f)$;
\item \ $(\widetilde f, \widetilde f)$ is loose but not by a small deformation.
\end{enumerate}
\end{thm}

This result settles a question raised in \cite[1.6]{K6}.  More
importantly, however, it relates the completeness question concerning
$\omega^\#$--invariants to fascinating recent work of D.\ Gon\c calves
and D.\ Randall. They produced many maps $\widetilde f \co S^m \to S^n$
which cannot be deformed away from themselves by small homotopies but
only via large deformations which use all the space available in $S^n$
(cf \ref{1.22} (6)).  For example all Whitehead products of the form
\begin{equation*} \
[\widetilde f] = [\iota_{4 k + 2}, \ \iota_{4 k + 2}] \ \in \ \pi_{8 k + 3} (S^{4 k + 2}), \ \ \ k = 1, 2, \dots \ \ ,
\end{equation*}
have this property (cf Gon\c calves and Randall
\cite{GR1}). Moreover the existence of such maps in dimension $(m, n)
= (4k - 2, 2 k), \ k > 4$, turns out to be equivalent to the Strong
Kervaire Invariant One Problem, ie, the existence of an element of
order 2 with Kervaire invariant one in the stable homotopy group
$\pi^S_{2 k - 2}$ (cf \cite{GR2}); examples of such maps exist in
dimension $2 k = 16, 32, 64$.

\section{Covering spaces and the pathcomponents of $E (f_1, f_2)$}\label{sec2}

Throughout this paper $M, N, \widetilde N$, and $Q$ will denote smooth connected manifolds (having the Hausdorff property and countable bases) without boundary, $M$ being compact; $f_1, f_2, f, \overline f, \ldots\co M \to N$ will be (continuous) maps.

In this section we study the set $\pi_0 (E (f_1, f_2))$ of
pathcomponents of $E (f_1, f_2)$ (cf \ref{1.6}) with the help of
coverings of $N$. We will need this when we discuss Nielsen
decompositions in \fullref{sec3}.

Consider the diagram
\begin{equation} \label{2.1}
\begin{CD}
@.                                             \widetilde N \\
@.                                              @VV{p}V \\
M             @>{ \ f_1, f_2 \ }>>                    N
\end{CD}
\end{equation}
where $p$ is a covering map. Pick base points $x_0 \in M, \ y_{0 i} := f_i (x_0) \in N, \ i = 1, 2$, and $\widetilde y_{0 2} \in p^{-1} (\{ y_{0 2}\} ) \subset \widetilde N$. For any path $\theta$ joining $y_{01}$ to $y_{02}$ in $N$ (ie $(x_0, \theta) \in E (f_1, f_2))$ define the subset
\begin{equation}\label{2.2}
\pi_\theta := \{ [ \theta^{- 1} (f_1 \scirc \sigma) \theta (f_2 \scirc \sigma)^{- 1}] \in \pi_1 (N, y_{0 2}) \mid [\sigma] \in \pi_1 (M, x_0) \ \}
\end{equation}
of $\pi_1 (N, y_{0 2})$, which consists of all homotopy classes of concatenated loops of the form
$$
\begin{CD}
y_{0 2} @>{\theta^{- 1}}>> y_{0 1} @>{f_1 \scirc \sigma}>> y_{0 1} @>{\theta}>> y_{0 2} @>{f_2 \scirc \sigma^{- 1}}>> y_{0 2} .
\end{CD}
$$

\begin{example}\label{2.3}\ 

\begin{enumerate}
\item[(i)] If $f_{i *} (\pi_1 (M)) = 0, \ i = 1, 2$ (eg, if $M$ is simply connected), then $\pi_{\theta} = \{ 0\}$;
\item[(ii)] if $\pi_1 (N)$ is abelian and $f_1 (x_0) = f_2 (x_0)$ then
$$
\pi_\theta \ = \ (f_{1 *} - f_{2 *}) (\pi_1 (M, x_0));
$$
\item[(iii)] (root case) if $f_2$ is constant, then $\pi_\theta$ is conjugated, via $[\theta]$, to $f_{1 *} (\pi_1 (M, x_0))$; whenever $\pi_1 (N)$ is not commutative, here (and in other cases) $\pi_\theta$ may depend strongly on $\theta$.
\end{enumerate}
\end{example}

Now consider the map
$$
\psi  \co\pr^{- 1} (\{ x_0\} ) = \{ (x_0, \theta) \in E (f_1, f_2)\} \longrightarrow p^{- 1} (\{ y_{0, 1} \} ) \ \subset \ \widetilde N
$$
defined on the fiber of the obvious projection $\pr \co E (f_1, f_2) \to
M$ (cf \ref{1.6}) as follows: let $\widetilde\theta \co I \to \widetilde
N$ be the unique (continuous) lifting of $\theta$ such that
$\widetilde\theta (1) = \widetilde y_{0 2}$ and put $\psi (x_0,
\theta) := \widetilde\theta (0)$. (The reason for \lq\lq lifting
backwards\rq\rq will be explained in remark \ref{4.5}.) There is also
the surjective map
$$
q\co \pr^{- 1}_* (\{ x_0 \} ) \longrightarrow \pi_0 (E (f_1, f_2))
$$
which assigns, to each element $(x_0, \theta)$, the corresponding pathcomponent of $E (f_1, f_2)$. When does $\psi$ factor through $q$?

\begin{prop}\label{2.4}
\ $\psi$ induces a welldefined bijection
$$
\Psi\co \pi_0 (E (f_1, f_2)) \longleftrightarrow \ p^{- 1} (\{ y_{0, 1} \} )
$$
if and only if
\begin{equation}\tag{$*$}\label{star}
\pi_\theta \ = \ p_* (\pi_1 (\widetilde N, \widetilde y_{0 2})) \ \ \text{for all} \ (x_0, \theta) \in E (f_1, f_2)\end{equation}
{\rm (}This condition holds for example if $M$ and $\widetilde N$ are simply connected{\rm )}.
\end{prop}

\begin{proof}
Let $(x_0, \theta'), (x_0, \theta)$ be elements of $E (f_1, f_2)$ based at $x_0$.

First assume that they can be joined by a path $(\sigma, \Theta)$ in $E (f_1, f_2)$. Then $\Theta$ yields homotopies from $\theta'$ to $(f_1 \scirc \sigma) \theta (f_2 \scirc \sigma)^{- 1}$ and from $\theta^{- 1} \theta'$ to $\theta^{- 1} (f_1 \scirc \sigma) \theta (f_2 \scirc \sigma)^{- 1}$ which leave the endpoints fixed. Clearly $\psi (\theta') = \psi (\theta)$ precisely if $\theta^{- 1} \theta'$ lifts to a closed loop in $\widetilde N$ starting and ending in $\widetilde y_{0 2}$. Thus $\psi$ induces a welldefined map on $\pi_0 (E (f_1, f_2))$ if and only if $\pi_\theta \subset p_* (\pi_1 (\widetilde N, \widetilde y_{0 2}))$ for all $(x_0, \theta) \in E (f_1, f_2)$. Surjectivity follows automatically since $\widetilde N$ is connected.

Next suppose only that $\psi (\theta') = \psi (\theta)$. If $\pi_\theta \supset p_* (\pi_1 (\widetilde N, \widetilde y_{0 2}))$, then $[\theta^{- 1} \theta'] = [\theta^{- 1} (f_1 \scirc \sigma) \theta (f_2 \scirc \sigma)^{-1}]$ for some $[\sigma] \in \pi_1 (M, x_0)$; this  yields a path in $E (f_1, f_2)$ joining $(x_0, \theta')$ to $(x_0, \theta)$. Thus condition \ref{star} in \ref{2.4} implies also the injectivity of $\Psi$.

On the other hand, given $(x_0, \theta) \in E (f_1, f_2)$ and $[\tau] \in p_* (\pi_1 (\widetilde N, \widetilde y_{0 2}))$, put $\theta' = \theta \tau$. If $\Psi$ is injective, then there is a path in $E (f_1, f_2)$ from $(x_0, \theta')$ to $(x_0, \theta)$, and we conclude again that $[\tau] = [\theta^{- 1} \theta']$ (and hence all of $p_* (\pi_1 (\widetilde N, \widetilde y_{0 2}))$) lies in $\pi_\theta$.

In the special case when $f_{i *} (\pi_1 (M)) = 0, \ i = 1, 2$, condition \ref{star} in \ref{2.4} holds if and only if $\widetilde N$ is simply connected (cf \ref{2.2} and \ref{2.3}).
\end{proof}

\section{Nielsen numbers and covering spaces}\label{sec3}

In this section we use liftings to covering spaces in order to give a new description of our Nielsen numbers. As a sample application we discuss maps into spherical space forms.

In the setting of diagram \ref{2.1} we assume that
\begin{enumerate}
\item[(i)] the group $G$ of covering transformations of the covering space $p \co \widetilde N \to N$ acts transitively on the fibers of $p$ (or, equivalently, $p_* (\pi_1 (\widetilde N))$ is a normal subgroup of $\pi_1 (N))$; and
\item[(ii)] there are liftings $\widetilde f_i \co M \to \widetilde N$ of $f_i$ (or, equivalently, $f_{i *} (\pi_1 (M) \subset p_* (\pi_1 (\widetilde N))),$ $ i = 1, 2$.
\end{enumerate}

Any choice of such maps $\widetilde f_i$ (satisfying $p \scirc \widetilde f_i = f_i$) determines a homeomorphism
\begin{equation}\label{3.1}
\lambda \co E (f_1, f_2) \ \ \longleftrightarrow \ \ \amalg_{g \in G} \ E (g \scirc \widetilde f_1, \widetilde f_2)
\end{equation}
(cf \ref{1.6}) defined by $\lambda (x, \theta) = (x, \widetilde\theta), \ (x, \theta) \in E (f_1, f_2)$, where $\widetilde\theta$ is the unique lifting of the path $\theta$ such that $\widetilde\theta (1) = \widetilde f_2 (x)$; here the disjoint components $E (g \scirc \widetilde f_1, \widetilde f_2)$ on the right hand side correspond to open subsets of $E (f_1, f_2)$.

$\lambda$, together with the obvious tangent isomorphism $T \widetilde N \cong p^* (TN)$, induces the map
\begin{equation}\label{3.2}
\lambda_* \co \Omega^\# (f_1, f_2) \ \longrightarrow  \ \prod_{g \in G} \Omega^\# (g \scirc \widetilde f_1, \widetilde f_2).
\end{equation}
(cf \ref{1.8}). Clearly
\begin{equation}\label{3.3}
\lambda_* (\omega^\# (f_1, f_2)) \ = \ (\omega^\# (g \scirc \widetilde f_1, \widetilde f_2))_{g \in G}
\end{equation}
(cf \ref{1.4}--\ref{1.8}) and $\lambda_* ([\phi]) = ([\phi])_{g \in
G}$. Possibly $\lambda_*$ is neither onto (eg when $G$ is infinite)
nor injective (since the disjointness requirements in the definition
of $\Omega^\# (f_1, f_2)$, concerning for example embedded bordisms, may not
be preserved). However, if $E (g \scirc \widetilde f_1, \widetilde
f_2)$ is pathconnected for all $g \in G$, ie, if the condition \ref{star}
in \ref{2.4} holds, then $\lambda_* (\omega^\# (f_1, f_2))$ keeps track of
strongly essential Nielsen components (compare \ref{1.9}). We conclude:

\begin{thm}\label{3.4}
Let $p \co \widetilde N \to N$ be a covering space such that the group $G$ of covering transformations acts transitively on the fibers and let $\widetilde f_i \co M \to \widetilde N$ be a lifting of $f_i, \ i = 1, 2.$ Assume that condition \ref{star}  {\rm (}of \fullref{2.4}{\rm )} holds.

Then $\# \pi_0 (E (f_1, f_2)) = \# G$ and the strong Nielsen number of $(f_1, f_2)$ {\rm (}cf \ref{1.9}{)} is given by
$$
N^\# (f_1, f_2) \ = \ \# \{ g \in G \mid \omega^\# (g \scirc \widetilde f_1, \widetilde f_2) \ \ne \ 0 \} .
$$
\end{thm}

It is easy to see that condition \ref{star} and the lifting condition are independent.   For example if $\widetilde f \co M \to \widetilde N$ and $p$ are nontrivial coverings and $\pi_1 (N)$ is abelian, then $ f_1, f_2 := p \scirc \widetilde f$ have liftings, but \ref{star} fails to hold (cf \ref{2.3} (ii) and \ref{2.4}).

On the other hand consider the case where $\pi_1 (\widetilde N) = 0$ and $f_1 (x_0) = f_2 (x_0) =: y_0$. Here $G$ acts transitively; moreover \ref{star} is satisfied if and only if $f_{1 *} = f_{2*}$ maps $\pi_1 (M, x_0)$ into the center of $\pi_1 (N, y_0)$; thus if for example $\pi_1 (N, x_0)$ is abelian and $f_1 = f_2$ is a non-universal covering map, then no liftings $\widetilde f_i$ exist, but condition \ref{star} holds.

In any case all the assumptions in \fullref{3.4} are satisfied whenever both $M$ and $\widetilde N$ are simply connected.

\begin{cor}[Spherical space forms]\label{3.5}
Given a free smooth action of a nontrivial finite group $G$ on $S^n, \ n \ge 1$, let $N = S^n / G$ be the quotient manifold. Consider maps $f_1, f_2 \co M \to S^n / G$ where {\rm (i)} $M = S^m, \ m \ge 2$; or {\rm (ii)} $M$ is simply connected, having dimension $m < 2 n - 2$.

Then $f_1$ and $f_2$ are homotopic whenever $N^\# (f_1, f_2) \ne  \# G$.

In particular, if $n$ is odd then
$$
N^\# (f_1, f_2) = MCC (f_1, f_2) = \begin{cases}
\ \# G   \qquad & \text{if} \ f_1 \not\sim f_2 \ ; \\
\ 0 & \text{if} \ f_1 \sim f_2  \ .
\end{cases}
$$
\end{cor}

\begin{proof}
\fullref{3.4} offers a new and more powerful approach when compared to the treatment in \cite{K6} (see statement 1.13 there and its proof (preceding 6.13) which works only if $f_1$ or $f_2$ is not coincidence producing, cf \ref{5.8} (iii) below).

Our claim holds trivially when $n = 1$ since then $f_1$ and $f_2$ are
nulhomotopic. Thus we may assume that $\pi_1 (S^n) = 0$. Suppose $N^\#
(f_1, f_2) < \# G$. Then there exist liftings $\widetilde f_1,
\widetilde f_2 \co M \to S^n$ such that $\omega^\# (\widetilde f_1,
\widetilde f_2) = 0$ (cf \ref{3.4}). Therefore $(\widetilde f_1, \widetilde
f_2)$ is loose; this follows from \ref{1.14} in case (i) and from \cite[1.2 (iii)]{K6},
 and \cite[1.10]{K3} in case (ii). After a homotopy
$\widetilde f_1, \widetilde f_2$ are coincidence free. Thus
$\widetilde f_1$ is homotopic to $A \scirc \widetilde f_2$, and so is
$g \scirc \widetilde f_2$ for every nontrivial element $g \in G$ (cf
Dold and Gon\c calves \cite[2.10]{DG}, or the beginning of Section 8
in \cite{K3}; here $A$ denotes the antipodal map). We conclude that
$$
f_1 \ = \ p \scirc \widetilde f_1 \ \sim \ p \scirc g \scirc \widetilde f_2 \ =  \ p \scirc \widetilde f_2 \ = \ f_2 .
$$

If $n$ is odd, $f_1$ can be pushed away from itself along a nowhere zero vector field. In view of \fullref{1.10} this completes the proof.
\end{proof}

\section{Fibrations and the $\omega^\#$--invariants}\label{sec4}

In this section we generalize $\lambda_*$ (cf \ref{3.2} and \ref{3.3}). We use this to explore compatibilities of our invariants with fibrations, eg, with the natural projections from spheres to complex or quaternionic projective spaces.

Consider commuting diagrams
\begin{equation}\label{4.1}
 \xymatrix{
    & Q \ar[d]^-{p} \\
    M \ar[ur]^-{\widetilde{f_i}} \ar[r]^-{f_i} & N\ ,
  }
\end{equation}
$i = 1, 2$, where $p$ is a smooth locally trivial fibration. We want to compare the $\omega^\#$--invariants of $(f_1, f_2)$ and of the pair $(\widetilde f_1, \widetilde f_2)$ of liftings.

Given a bordism class
$$
c \ = \ [ \ C \subset M, \ \widetilde g \co C \to E (f_1, f_2), 
   \ \overline g^\#] \ \in \ \Omega^\# (f_1, f_2)
$$
(compare \ref{1.4}--\ref{1.8}), let $H \co C \times I \to N$ be the
homotopy from $f_1 | C$ to $f_2 | C$ which is adjoint to $\widetilde
g$, and let $\widetilde H \co C \times I \to Q$ be a lifting of $H$
which {\em ends} at $\widetilde H (\ , 1) = \widetilde f_2 | C$. Then
$\widetilde f_1 | C$ and $\widetilde H (\ , 0)$ determine sections of
the pulled back fibration $(f_1 | C)^* (p)$. Consider their transverse
intersection locus $C' \subset C \subset M$. The normal bundle $\nu
(C', C)$ is canonically isomorphic to $(\widetilde f_1 | C')^* (TF)$,
the pullback of the tangent bundle along the fibers of $p$; together
with $\overline g^\# \co \nu (C, M) \cong f_1^* (TN) | C$ this yields a
vector bundle isomorphism
\begin{equation}\label{4.2}
\overset{\simeq}{g}^\# \co  \nu (C', M) \cong \widetilde f^*_1 (TQ) | C' \
\end{equation}
(compare \ref{1.5}).
Moreover the adjoint of $\widetilde H | (C' \times I)$ determines a section $\overset{\approx}{g}$ over $C'$ of the fibration $E (\widetilde f_1, \widetilde f_2) \to M$ (compare \ref{1.6}). We obtain the bordism class
\begin{equation}\label{4.3}
\lambda_{e *} ([c]) \ = \ [ \ C' \subset M, \ \overset{\approx}{g}, \ \overset{\simeq}{g}^\# ] \ \ \in \ \ \Omega^\# (\widetilde f_1, \widetilde f_2) \ .
\end{equation}

\begin{prop}\label{4.4}
This construction yields a welldefined map
$$
\lambda_{e *} \co \Omega^\# (f_1, f_2) \longrightarrow \Omega^\# (\widetilde f_1, \widetilde f_2) 
$$
which takes $\omega^\# (f_1, f_2)$ to $\omega^\# (\widetilde f_1, \widetilde f_2)$ and $[\phi]$ to $[\phi]$. In particular, if $\omega^\# (f_1, f_2)$ is trivial then so is $\omega^\# (\widetilde f_1, \widetilde f_2)$.

In the special case where $p$ is a covering map as in \ref{3.1} and \ref{3.2} $\lambda_{e *}$ is the component map of $\lambda_*$  corresponding to the unit $e$ of the group $G$ of covering transformations.
\end{prop}

The proof is fairly straight forward.

\begin{remark}\label{4.5}
In the definitions of $\psi$ (cf \ref{2.4}), $\lambda$ (cf \ref{3.1}), and $\lambda_{e *}$ (cf \ref{4.3}) we lifted paths and homotopies \lq\lq backwards\rq\rq, starting at the end. This makes our constructions more easily compatible with our convention to describe normal bundles such as $\nu (C, M)$ in terms of $f^*_1 (TN)$ (and not of $f^*_2 (TN)$, cf \ref{1.5}) and helps us to avoid reframings as the ones necessitated for example in the definitions of $\alpha$ and $\beta$ in \cite[(48)--(50)]{K6}.
\end{remark}

The following sample application of \fullref{4.4} will be needed in the proof of 
\fullref{1.17}.

\begin{cor}\label{4.6}
Let $p \co S^{d n' + d - 1} \to \mathbb K P (n'), \ n' \ge 1$, be the canonical fibration where $\mathbb K = \mathbb C$ or $\mathbb H$ {\rm(}with real dimension $d = 2$ or $4)$. Assume that $M = S^m$ or that $M$ has dimension $m < 2 d n' - 2$. If the maps $f_i \co M \to \mathbb K P (n')$ allow genuine liftings $\widetilde f_i$ to $S^{d n' + d - 1}$ (ie $p \scirc \widetilde f_i = f_i), \ i = 1, 2$, and if $\omega^\# (f_1, f_2) = 0$, then $f_1$ and $f_2$ are homotopic.
\end{cor}

\begin{proof}
If $\omega^\# (\widetilde f_1, \widetilde f_2) = 0$ (in view of \ref{4.4})
then the pair $(\widetilde f_1, \widetilde f_2)$ is loose and
$\widetilde f_1 \sim A \scirc \widetilde f_2$ (cf \cite[1.12,
1.14]{K6}, \cite[1.2 (iii)]{K3}, 1.10 and \cite[2.10]{DG}). Therefore
$f_1 \sim p \scirc A \scirc \widetilde f_2 = p \scirc \widetilde f_2 =
f_2$.
\end{proof}

\section{Removing selfcoincidences}\label{sec5}

In view of results such as Corollaries \ref{3.5} and \ref{4.6} a closer look at
the selfcoincidence setting seems to be in order. Here we are dealing
essentially with pairs of the form $(f, f)$. Only the pathcomponent of
$E (f, f)$ which contains all elements of the form $(x, \
\text{constant path})$ (cf \ref{1.6}) can possibly be strongly essential (cf
\ref{1.9}). As a consequence our invariants $\omega^\#$ and $N^\#$ turn out
to be not as powerful here as in the root case (where $f_1$ or $f_2$
is constant, compare eg \cite[6.5 b) (iv)]{K6}). Furthermore the
second part of the first question in \ref{1.13} is easily answered.

\begin{prop}\label{5.1}
Let $f_1, f_2 \co M \to N$ be homotopic. Then  $MCC (f_1, f_2) \le 1$; moreover $\omega^\# (f_1, f_2)$ is trivial if and only if $N^\# (f_1, f_2) = 0$. In the special case where $M$ is a sphere $MC (f_1, f_2)$ equals $MCC (f_1, f_2)$.
\end{prop}

\begin{proof} \
The coincidence space $C (f_1, f_1) = M$ is pathconnected by assumption. If $(f_1, f_2)$ is a generic pair of sufficiently close maps then $C_A $ (cf \ref{1.9} (i)) is empty for all but possibly one pathcomponent $A$ of $E (f_1, f_2)$; the corresponding partial bordism class, ie $\omega^\# (f_1, f_2)$, vanishes precisely when $N^\# (f_1, f_2) = 0$. If  $f_1$ and $f_2$ are related by an arbitrary (possibly \lq\lq large\rq\rq ) homotopy, apply \cite[2.1]{K6}.

If $M = S^m$ the vector bundle $f^*_1 (TN)$ allows a section with at most one zero. This yields a map $f'_2$ such that $\# C (f_1, f'_2) \le 1$ and $f'_2 \sim f_2$.
\end{proof}

\medskip
Let $\xi \subset TN$ be a sub--vectorbundle and let $\xi'$ denote the image of its total space under the composed diffeomorphism
\begin{equation}\label{5.2}
\begin{CD} TN  @<{\cong}<{p_{1*}}< \nu (\Delta, N \times N) \ \cong \ U \ \ \ \subset \ \ \ N \times N \end{CD}
\end{equation}

where $U$ is a tubular neighbourhood of the diagonal $\Delta$ in \ $N \times N$ and $p_{1 *}$ denotes the vector bundle isomorphism induced by the first projection $p_1 \co N \times N \to N$.

\begin{definition}\label{5.3}
Given a map $f \co M \to N$, we say $(f, f)$ is {\em loose by a small}\break
$\xi$--{\em deformation} if for every metric on $N$ and for every
$\epsilon > 0$ there exists an $\epsilon$--approximation $\overline f
\co M \to N$ of $f$ such that $(f, \overline f) (M)\ \subset \ \xi' -
\Delta$. If this holds for $\xi = TN$ we simply say that $(f, f)$ is
{\em loose by a small deformation} (compare \cite{DG} or \cite[Section
1]{GR2}).
\end{definition}

A homotopy lifting argument shows that $(f, f)$ is loose by a small $\xi$--deformation precisely if the pulled back vector bundle $f^* (\xi)$ has a nowhere vanishing section over $M$ (compare \cite{DG}).

If $M = S^m$ we can approach such phenomena using standard tools of homotopy theory.

Assume $m, n \ge 2$ and consider the commuting diagram
\medskip
\begin{equation}\label{5.4}
  \xymatrix{
    && \pi_m(S^n) \\
    \dots \pi_m (STN) \ar[r] \ar[d] &
    \pi_m (N)  \ar[r]^-{\partial} \ar@{=}[d] \ar@{-->}[ur]^-{\underline\omega^\#}&
    \pi_{m -1} (S^{n -1}) \ar[r]^-{\incl_*} \ar[d]^-{j_*} \ar[u]_-{E}&
    \pi_{m -1} (STN) \ar[d] \\
    \dots \pi_m \left(\widetilde C_2 (N)\right)  \ar[r] &
    \pi_m (N) \ar[r]^-{\partial '} &
    \pi_{m-1} \left(N - \{ x_0\}\right)  \ar[r]^-{\incl'_*} &
    \pi_{m -1} \left(\widetilde C_2 (N)\right)
  }
\end{equation}

where the fibrations of the space $STN$  of unit tangent vectors and of the configuration space
\begin{equation}\label{5.5}
\widetilde C_2 (N) \ = \ \{ (y_1, y_2) \in  N\times N\mid y_1 \ne y_2\} \ = \ N \times N - \Delta
\end{equation}
over $N$ yield the exact horizontal sequences;  $\incl, \ \incl'$ \ (and $E$, resp.) denote fiber inclusions (and the Freudenthal suspension, resp.). The downward pointing vertical arrows are induced by a diffeomorphism as described in \ref{5.2}; in particular, $j$ denotes the inclusion of the boundary sphere of an $n$--ball in $N$ around the basepoint $x_0$.

Given a homotopy class $[f] \in \pi_m (N)$, how does its image under the boundary operator $\partial$ in \ref{5.4} compare to our looseness obstruction $\omega^\# (f, f)$? Recall that the path space approach yields no added information in the selfcoincidence setting; hence $\omega^\# (f, f)$ is just as strong as the invariant
\begin{equation}\label{5.6}
\underline{\omega}^\# (f, f) \ := \ \pr_* (\omega^\# (f, f)) \ \in \pi_m (S^n)
\end{equation}
defined by the framed bordism class of the zero set $Z$ of a generic
section $s$ of $f^* (TN)$ (cf \cite[(43)--(45) and section 6]{K6}). We
may assume that $Z$ lies in the interior of a small ball $B \subset
S^m$ (over which $f^* (TN)$ is trivialized) and that $s$ maps the
boundary $\partial B \cong S^{m -1}$ to unit vectors (and thus defines
a map $s| \co S^{m -1} \to S^{n -1})$. Then $\partial ([f]) = [s|\ ]$
and $Z$ is framed bordant to $s|^{-1} (\{ *\} )$, $*$ a regular value
of $s|$. We conclude:
\begin{prop}\label{5.7}
For all $[f] \in \pi_m (N)$, we have
\begin{equation*}
\underline{\omega}^\# (f, f) \ = \ \pm E  \scirc  \partial ([f]) .
\end{equation*} 
\end{prop}

Given $[f] \in \pi_m (N)$, consider the following conditions:

\begin{numberonly}\label{5.8}\ 

\begin{enumerate}
\item[(i)] $(f, f)$ is loose by a small deformation; equivalently, $\partial ([f]) = 0$;
\item[(ii)] $(f, f)$ is loose (by any deformation);
\item[(iii)] there exists any map $\overline f \co S^m \to N$ such that $(f, \overline f)$ is loose (we say that $f$ is {\em not coincidence producing}, cf 
Brown and Schirmer \cite{BS}); equivalently, $\partial' ([f]) = 0$;
\item[(iii$'$)] $\omega^\# (f, f) = 0$.
\end{enumerate}
\end{numberonly}

Clearly (i) implies (ii). In turn, (ii) implies both (iii) and (iii$'$).

\begin{prop}\label{5.9}
The conditions {\rm (i)} and {\rm (iii)} are equivalent for all $[f] \in \pi_m (N)$ if and only if the homomorphism
\begin{equation*}
(j_*, \incl_*) \co \ \pi_{m -1} (S^{n -1}) \ \longrightarrow \ \pi_{m -1} (N - \{ x_0\} ) \oplus \pi_{m -1} (STN)
\end{equation*}
is injective (where $\incl$ denotes the fiber inclusion).
\end{prop}

Indeed, this is the precise condition for the kernels of $\partial$ and $\partial' = j_* \scirc \partial$ in \ref{5.4} to agree. In view of \ref{5.6} and \ref{5.7} we obtain similarly

\begin{prop}\label{5.10}
The conditions {\rm (i)} and {\rm (iii$'$)} are equivalent for all $[f] \in \pi_m (N)$ if and only if the homomorphism
\begin{equation*}
(E, \incl_*) \co \ \pi_{m -1} (S^{n -1}) \ \longrightarrow \ \pi_m (S^n) \oplus \pi_{m -1} (STN)
\end{equation*}
is injective.
\end{prop}

In particular, if $j_*$ is injective on $\im \partial = \ker \incl_*$ but $E$ is not, then there is a map $f \co S^m \to N$ such that $\omega^\# (f, f) = 0$ but $(f, f)$ is not loose by any deformation. In the next section we will study such examples systematically in case $N$ is a projective space. They show that even the nonstabilized \lq\lq desuspended\rq\rq invariant $\omega^\#$ (compare \cite{K3}) can be (at least one desuspension) short of yielding complete looseness criteria.

\section{Maps into projective spaces}\label{sec6}

In this section we apply the previous discussion to the case $M = S^m,
\ N = \mathbb K P (n')$ (where $\mathbb K = \mathbb R, \ \ \mathbb C  \ \text{or} \ \mathbb H$), exploiting natural compatibilities with the canonical fiber map
\begin{equation}\label{6.1}
p  \co  S^{n + d - 1} \longrightarrow \mathbb K P (n') \ .
\end{equation}
This leads to a proof of Theorems \ref{1.17}, \ref{1.21}, and \ref{1.22}.

First we show that the inequalities of \fullref{1.10} hold here in full generality, without an exception at the dimension $n = 2$ which is so critical in classical fixed point and coincidence theory (cf Jiang \cite{Ji1,Ji2} and also \cite{BGZ1,BGZ2}, and \cite{GKZ}).

\begin{prop}\label{6.2}
For all maps $f_1, f_2 \co S^m \to \mathbb K P (n')$ {\rm (}where $m, n' \ge 1$ and $\mathbb K = \mathbb R, \ \mathbb C \ \text{or} \ \mathbb H)$ we have:
\begin{enumerate}
\item[{\rm (i)}] $MCC (f_1, f_2) \ \le \ \# \pi_0 (E (f_1, f_2))$; \ and
\item[{\rm (ii)}] $MC (f_1, f_2) \ \le \ \# \pi_0 (E (f_1, f_2))$ or $MC (f_1, f_2) = \infty$.
\end{enumerate}
\end{prop}

\begin{proof}
In view of \ref{1.10}, \ref{1.14}, \ref{1.15}, and Jezierski \cite[4.0]{Je}, 4.0 (see
also \cite[6.14]{K6}) we need to consider only the case when $m >
2$, $N = \mathbb R P (2)$. Then $\pi_{m -1} (S^{n -1}) = 0$ and hence
$(f_2, f_2)$ is loose (cf \ref{5.4} and \ref{5.8}).

Thus $MCC (f_1, f_2) = MCC (f, *)$ where $[f] = [f_1] - [f_2]$ (by
\cite[6.2]{K6}). Given liftings $\widetilde f, \widetilde * \co S^m \to
S^2$ of $f$ and of the constant map, $C (\widetilde f, \widetilde *) =
\widetilde f^{- 1} (\{ *\} )$ is generically a framed submanifold
which we may make connected by a suitable surgery. Taking also the
inverse image of a nearby (and, after an isotopy, antipodal point) we
see that $f^{- 1} (\{ *\})$ consists of two \lq\lq parallel\rq\rq\
connected components. Hence $MCC (f, *) \ \le 2 \ = \ \# \pi_0 (E
(f_1, f_2))$ (cf \cite[2.1]{K3}).
\end{proof}

For the remainder of this section we assume that $m, n' \ge
2$. Comparing the exact homotopy sequences of the fibrations $p$ (cf
\ref{6.1}) and $p| \co S^{n -1} \to \mathbb K P (n' - 1)$ we see that $\pi_m
(\mathbb K P (n')) $ is the direct sum of $p_* (\pi_m (S^{n + d -
  1}))$ and of the image of  $\pi_m (\mathbb K P (n' -1))$; moreover
$p_*$ is injective. Thus
\begin{equation}\label{6.3}
[f_i] \ \ = \ \ [p \scirc \widetilde f_i] \ + \ [f''_i] \ \
\end{equation}
where $(f''_i, *)$ is loose (for $* \not\in \mathbb K P (n' -1)), \ i = 1, 2$. Clearly the pairs $(f_1, f_2)$ and $(f'_1, f'_2) := (p \scirc \widetilde f_1, p \scirc \widetilde f_2)$ have identical $\omega^\#$--invariants as well as Nielsen and minimum numbers (cf \cite[6.1 and 6.2]{K6}). Moreover since $f''_i$ factors through $\mathbb K P (n' -1)$ the corresponding pullback of the sphere bundle $S T (\mathbb K P (n'))$ allows a section; thus $\partial ([f''_i]) = 0$ and
\begin{equation}\label{6.4}
\partial ([f_i]) \ \ = \ \ \partial ([p \scirc \widetilde f_i]) \ .
\end{equation}
This reduces the proof of \fullref{1.17} to the case where $f_i = p \scirc \widetilde f_i = f'_i,\ i = 1, 2$.

Now apply the discussion of the diagram \ref{5.4} to the case $N = \mathbb K P (n')$. Clearly $j_*$ is injective here since $j$ agrees with $p| S^{n -1}$ up to the homotopy equivalence $\mathbb K P (n') - \{ x_0\} \sim \mathbb K P (n' -1)$. Thus the three looseness conditions (i), (ii), and (iii) in \ref{5.8} are equivalent for all $[f] \in \pi_m (\mathbb K P (n')$ (cf \ref{5.9}). They hold precisely if $\partial ([f]) = 0$ and, in case $[f] = p_* ([\widetilde f])$, precisely if $\partial_{\mathbb K} ([\widetilde f]) = 0$. Indeed, the exact horizontal sequences in the diagrams \ref{1.16} and \ref{5.4} are closely related via maps induced by $p$; in particular,
\begin{equation}\label{6.5}
\partial_{\mathbb K} \ \ \ = \ \ \ \partial \ \scirc \ p_* \ \ \ .
\end{equation}
Here we observe that the vector bundle $p^* (T \mathbb K P (n'))$ is isomorphic to the orthogonal complement
\begin{equation}\label{6.6}
\xi_{\mathbb K} \ := \ \{ (x, v) \ \in \ S^{n + d - 1} \times \mathbb
K^{n' + 1} \mid v \perp \mathbb K \cdot x \ \} \ \ \ \ \ \ \subset \
\ TS^{n + d - 1}
\end{equation}
of the $\mathbb K$--line bundle over $S^{n + d - 1}$ which is spanned by the locus vectors;  eg $\xi_{\mathbb R} = T S^n$. The corresponding space $S (\xi_{\mathbb K})$ of unit vectors is the Stiefel manifold $V_{n' + 1, 2} (\mathbb K)$.

In view of \fullref{5.1} the following result implies \fullref{1.22}.

\begin{thm}\label{6.7}
Assume $m, n' \ge 2$ and $\mathbb K = \mathbb R, \ \mathbb C$ or $\mathbb H$. Given $[f] \in \pi_m (\mathbb K P (n'))$, let $[\widetilde f] \in \pi_m (S^{n + d - 1})$ be the corresponding component in the decomposition  {\rm \ref{6.3}}. Then the following conditions are equivalent:
\begin{enumerate}
\item[{\rm (i)}] $(f, f)$ is loose by a small deformation;
\item[{\rm (ii)}] $(f, f)$ is loose by any deformation;
\item[{\rm (iii)}] $f$ is not coincidence producing {\rm (}cf {\rm \ref{5.8}, iii}{\rm )};
\item[{\rm (iv)}] $(\widetilde f, \widetilde f)$ is loose by a small $\xi_{\mathbb K}$--deformation {\rm (}cf {\rm \ref{5.3})};
\item[{\rm (v)}] $\partial_{\mathbb K} ([\widetilde f]) = 0$.
\end{enumerate}
{\rm (}Note that \lq\lq loose by a small $\xi_{\mathbb R}$--deformation\rq\rq\ is the same as \lq\lq loose by a small deformation\rq\rq {\rm )}.

Furthermore, consider also the following {\rm (}possibly weaker{\rm )} conditions:
\begin{enumerate}
\item[{\rm (iii$'$)}] $\omega^\# (f, f) = 0$ ;
\item[{\rm (iii$''$)}] $\omega^\# (\widetilde f, \widetilde f) = 0$ ;
\item[{\rm (iv$'$)}] $(\widetilde f, \widetilde f)$ is loose; \ \ \ \ \ \ and
\item[{\rm (v$'$)}] $E (\partial_{\mathbb K} ([\widetilde f])) = 0$ \ .
\end{enumerate}

If $\mathbb K = \mathbb R$, then all four conditions  {\rm (iii$'$) -- (v$'$)} are equivalent. If $\mathbb K = \mathbb C$ or $\mathbb H$, we still have the following implications:
\begin{equation*}
(\text{\rm v$'$})  \ \ \Longleftrightarrow \ \ (\text{\rm iii$'$}) \ \ \Longrightarrow \ \ (\text{\rm iii$''$}) \ \ \Longleftrightarrow \ \ (\text{\rm iv$'$}) \ .
\end{equation*}
\end{thm}

\begin{proof} \
The first claim follows from the previous discussion and is based essentially on \ref{5.9}, \ref{1.16}, and the remark following \ref{5.3}.

If $\mathbb K \, = \, \mathbb R$ then $\widetilde f$ is a genuine lifting of $f$ to a double cover. Clearly $\underline\omega^\# (f, f) = \underline\omega^\# (\widetilde f, \widetilde f)$ (cf \ref{5.6}). But these invariants contain just as much information as $\omega^\# (f, f)$ and $\omega^\# (\widetilde f, \widetilde f)$. In view of \ref{1.14}, \ref{4.4}, \ref{5.7}, and \ref{6.5} this completes the proof.
\end{proof}

\begin{proof}[Proof of Theorems \ref{1.17} and \ref{1.21}.] 
In view of the discussion of \ref{6.3} we may assume that $f_i = f'_i  = p \scirc \widetilde f_i, \ i = 1, 2$. Also note that
\begin{equation*}
\# \pi_0 (E (f_1, f_2)) = \ \# \pi_1 (\mathbb K P (n'))  = \begin{cases}\ 2 \qquad & \text{if} \ \mathbb K \, = \, \mathbb R \ \ \ , \\
\ 1  \qquad & \text{if} \ \mathbb K \, = \, \mathbb C \  {\text{or}}  \ \mathbb H \ \ , \\  \end{cases}
\end{equation*}
since $\pi_1 (S^m) = 0$ (cf \cite[2.1]{K3}).

First consider the case where $N^\# (f_1, f_2) < \# \pi_0 (E (f_1, f_2))$. Then $f_1 \sim f_2$ and $MCC (f_1, f_2) = MC (f_1, f_2) \le 1$ (cf \ref{3.5}, \ref{4.6}, and \ref{5.1}).  According to \fullref{1.10} (and \ref{5.1};  \ref{3.4}, \ref{1.14}) we can distinguish (and characterize) the three special cases:

\begin{enumerate}
\item $N^\# (f_1, f_2) = MCC (f_1, f_2) = 0$ \  (ie $f_1 \sim f_2$ and $(f_2, f_2)$ is loose);
\item $N^\# (f_1, f_2) < MCC (f_1, f_2) = 1$ \  (or, equivalently, $f_1 \sim f_2$ and $\omega^\# (f_2, f_2) = 0$ but $(f_2, f_2)$ is not loose);
\item $N^\# (f_1, f_2) = MCC (f_1, f_2) = 1$ \  (or, equivalently,  $\mathbb K = \mathbb R, f_1 \sim f_2$ and $\widetilde f_2$ is not homotopic to $A \scirc \widetilde f_2$).
\end{enumerate}

In view of \fullref{6.7} (and of the injectivity of $p_*$) these are just the first three cases in \fullref{1.18}.

In the remaining cases $N^\# (f_1, f_2)$ is equal to $\# \pi_0 (E
(f_1, f_2))$ and hence also to $MCC (f_1, f_2)$ and to $MC (f_1, f_2)$
(unless $MC (f_1, f_2) = \infty$, cf \ref{1.10} and \ref{6.2}). If $\mathbb K =
\mathbb C$ or $\mathbb H$ and $f_1, f_2$ have only finitely many
coincidence points we may deform $\widetilde f_1$ away from
$\widetilde f_2$ along the fibers of $p$; hence $\widetilde f_1 \sim A
\scirc \widetilde f_2, f_1 \sim f_2$ and $MC (f_1, f_2) \le 1$ (cf \ref{5.1}
and \cite[2.10]{DG}). Finally assume $\mathbb K = \mathbb
R$. According to \ref{5.1} and \cite[6.2]{K6}, $MC (f_1, f_2) < \infty$ if
and only if $MC (f_1 - f_2, *)$ (and hence $MC (\widetilde f_1 -
\widetilde f_2, \widetilde *))$ is finite. Then $\widetilde f_1 -
\widetilde f_2$ is a suspended map (cf \cite[1.12 or
6.10]{K6}). Conversely any suspended map (in the unreduced sense) maps
only the two suspension poles in $S^m$ to the suspension poles in
$S^n$, ie, to the fiber of the corresponding point in $\mathbb R P
(n)$. This establishes \fullref{1.18}.
\end{proof}

\bibliographystyle{gtart}
\bibliography{link}

\end{document}